\documentclass{article}
\usepackage{arxiv}
\usepackage[utf8]{inputenc} 
\usepackage[T1]{fontenc}    
\usepackage{hyperref}       
\usepackage{url}            
\usepackage{booktabs}       
\usepackage{amsmath}
\usepackage{amsthm}
\usepackage{amssymb}
\usepackage{mathtools}
\usepackage{amsfonts}       
\usepackage{nicefrac}       
\usepackage{microtype}      
\usepackage{graphicx}
\usepackage{doi}
\usepackage[textsize=tiny]{todonotes}

\usepackage{bbm}
\usepackage{siunitx}
\usepackage{caption}
\usepackage{mathrsfs}  
\usepackage{subcaption}
\usepackage{xcolor}
\usepackage{tcolorbox}
\usepackage[ruled]{algorithm2e}

\usepackage[style=alphabetic,sorting=nyt,sortcites=true,maxbibnames=9,autopunct=true,babel=hyphen,hyperref=true,abbreviate=false,backref=true,backend=biber]{biblatex}
\newtheorem{theorem}{Theorem}
\newtheorem{definition}{Definition}

\newtheorem{corollary}{Corollary}
\newtheorem{remark}{Remark}
\newtheorem{lemma}{Lemma}

\DeclareMathOperator*{\argmin}{arg\,min}

\newgeometry{vmargin={25mm}, hmargin={35mm,40mm}}

\title{Dynamical low-rank approximations of solutions to the Hamilton-Jacobi-Bellman equation}

\author{
	\textbf{Martin Eigel} \\
	Weierstrass Institute for \\ Applied Analysis  and Stochastics \\
	Berlin, Germany \\
	\texttt{eigel@wias-berlin.de}
	\And
	\textbf{Reinhold Schneider} \\
	Department of Mathematics \\
	Technical University Berlin \\
	Berlin, Germany \\
	\texttt{schneidr@math.tu-berlin.de}
	\And
	\textbf{David Sommer} \\
	Weierstrass Institute for \\ Applied Analysis  and Stochastics \\
	Berlin, Germany \\
	\texttt{sommer@wias-berlin.de} 
}

\renewcommand{\shorttitle}{DLR approximations of solutions to the HJB equation}

\begin{document}
\maketitle

\begin{abstract}
We present a novel method to approximate optimal feedback laws for nonlinear optimal control based on low-rank tensor train (TT) decompositions.
The approach is based on the Dirac-Frenkel variational principle with the modification that the optimisation uses an empirical risk.
Compared to current state-of-the-art TT methods, our approach exhibits a greatly reduced computational burden while achieving comparable results.
A rigorous description of the numerical scheme and demonstrations of its performance are provided.
\end{abstract}


\section{Introduction}
\label{sec:introduction}

Feedback control is ubiquitous and indispensable in real dynamical systems.
Since the controlled system can in general not be expected to follow model predictions exactly, system trajectories will eventually leave the forecasted path, meaning that any preplanned series of controls (albeit an optimal one) is based on wrong assumptions and therefore not only suboptimal, but potentially dangerous.
As an illustration, one might think of an astronaut who calculated an optimal course to land on the moon but then does not modify the forecasted actuation values of their rocket-drive when atmospheric effects steer them off said course, which will leave them drifting to outer space.
It is therefore vital to deploy controls based on \textit{current} state feedback, where \textit{current} means as frequently as possible in practical applications.

However, The problem of computing an optimal feedback control law for nonlinear optimal control problems is notoriously difficult.
This is because the synthesis of such a feedback law requires solving the Hamilton-Jacobi-Bellman (HJB) equation, which is a nonlinear parabolic partial differential equation (PDE) 
of generally high dimension $d \gg 1$ \cite{BC97}.
Classical schemes to solve the HJB equation such as Galerkin-schemes in linear ansatz spaces suffer from the \textit{curse of dimensionality} \cite{KK18}, i.e. an exponential complexity growth.
In practice, this means that the computation of a solution is often infeasibly slow if it can be discretised and stored at all.
Another severe obstacle can be the low regularity of viscosity solutions, cf~\cite{BD97}.
In this paper, our focus lies on the alleviation of the curse of dimensionality in order to enable the numerical treatment of high-dimensional control problems.
We hence only consider problems where the lack of regularity is not present or not pronounced enough to prevent a sufficiently accurate approximation.

The relevance of efficient numerical methods can be seen by the fact that true feedback control methods - that is: methods solving the HJB equation - are rarely used in practice due to the necessary computational effort.
Control problems arising e.g. in mechanical engineering often require new planning of controls to be computed within seconds.
There hence is a tight upper limit on the time budget available for generating new controls.
Therefore, most engineers deploy variations of Model Predictive Control (MPC) where open-loop controls are computed in such rapid succession that they effectively \textquotedblleft close the loop\textquotedblright ~\cite{camacho2013model}.
This is a conservative approach since the feedback property of the resulting controller is obtained solely by means of the measurements at the discrete planning steps.
In between two state measurements, the controller is not in feedback form.

In this work, we present a novel method to tackle nonlinear optimal control problems that $(1)$ yields a true feedback controller and $(2)$ has greatly reduced computational cost compared to current state-of-the-art methods.
Our method is based on policy iteration, linearising the HJB equation (which is then sometimes called the generalised Hamilton-Jacobi-Bellman or GHJB equation) 
and a modification of the Dirac-Frenkel variational principle.
This then allows the computation of approximate solutions on a specified function manifold, for which we choose the set of multivariate polynomials with a fixed tensor train (TT) rank.

Tree based tensor networks and tensor trains in particular have already been used for successful approximations of the value function in various works, see e.g. \cite{OS21,DKK19,FS20}. 
These recent results are summarized in the PhD thesis of Leon Sallandt \cite{Sallandt_2021}, which is still being finalised as this paper is written.
There, the approach is based on a Lagrangian (or dynamic programming) perspective by computing the value functions for several initial states and learning the global function from these values by regression using a multi-polynomial TT model.
With appropriate modifications, this approach can already be combined with regression techniques performed e.g. by machine learning methods, in particular artificial neural networks (NN). 
In the present paper, we follow a different approach, exploiting the Riemannian structure of the TT manifold ~\cite{Sch10,Steinlechner16} by an empirical version of the Dirac-Frenkel principle.

The solution obtained by the abstract Dirac-Frenkel principle can be shown to be quasi-optimal in some time interval $[0,T_\text{DF}]$ but deteriorates from the best low-rank approximation after a certain time \cite{LRS13}.
We expect a similar behaviour in our case which may restrict the time interval in practice. Combining both approaches - abstract and empirical - is an open research question, which we aim to address in future work.
Similarly, we defer the stochastic control case to a forthcoming paper, confining ourselves to deterministic control in this paper. We conjecture that the present approach is even more advantageous in the stochastic case.


The rest of the paper is organised as follows:
In Section~\ref{sec:rel_lit} we provide a short overview of the related literature, specifically the current state-of-the-art of tensor based methods to solve the HJB equation.
Section~\ref{sec:oc} introduces the finite horizon optimal control problem in feedback form, which the rest of this work revolves around.
In Section~\ref{sec:tt} the tensor train format, the corresponding manifold and the representation of the tangent space are introduced.
These are needed to formulate the Dirac-Frenkel variational principle, which is introduced in its abstract form in Section~\ref{sec:dlra}.
In Section~\ref{sec:pi}, we combine the concepts of Sections~\ref{sec:tt} and~\ref{sec:dlra} to develop our proposed DLRA method for approximately solving the HJB equation.
Numerical results that illustrate the practical performance are presented and discussed in Section~\ref{sec:num}.
Finally, we close in Section~\ref{sec:close} with an outlook on future work.

\section{Related work}
\label{sec:rel_lit}

The Bellman equation governing the value function of an optimal control (OC) problem was introduced as early as 1957 by Richard Bellman~\cite{Bell57}.
Since then, numerous sophisticated methods have been introduced to approximate solutions, mostly based on the principle of dynamical programming, see e.g.~\cite{Bert05} for a broad introduction to the subject.
The alternative approach, which we follow in this work, is to consider the infinitesimal version of the Bellman equation, namely the Hamilton-Jacobi-Bellman equation~\cite{BC97}, which is a nonlinear parabolic PDE.
In both cases, many methods rely on a fixed point iteration of the equation, which in the OC and Reinforcement Learning (RL) literature is called \textit{policy iteration}~\cite{How60}.
Alternatives are domain splitting algorithms~\cite{FLA94}, semi-Lagrangian methods~\cite{F87,FK14,TAK17}, data-based methods using Neural Networks~\cite{L14}, variational iterative methods~\cite{KDK13}, actor-critic methods~\cite{ZH21}, tree-based methods~\cite{AS19} and tropical methods~\cite{A08,A18}.

For a fixed starting value, an optimal control can be obtained by open-loop approaches such as Pontryagin's maximum principle~\cite{Pont61,Pont87}.
In this way, the value function can be evaluated pointwise by simply adding up the cost of that optimal control.
Controls of this type have been used to find the value function e.g. in~\cite{KW17,NTG19,BKK21,OS21}.
In this work, we use optimal open-loop controls as benchmarks to which we compare the feedback controller computed by our method.

Since any solution method for the HJB equation has to deal with the curse of dimensionality, some form of model order reduction has to take place in practical applications.
Possible function approximations can be obtained by using neural networks~\cite{D19,NR21,Ito21} or sparse polynomials~\cite{BKK21}.
In this work we use the TT format introduced to the mathematical community by Oseledets~\cite{O09,O11} for multivariate polynomials. A striking recent example of the power of the low rank TT structure for function approximation can be found in \cite{RSN21}, in which the authors use TTs with polynomial basis functions to outperform state-of-the-art NNs on the solution of parabolic PDEs by orders of magnitude, while requiring lower computational time.
For further details on TTs and more general hierarchical tensor networks, we refer the reader to the survey articles~\cite{HS14,BUS16} and the standard textbooks~\cite{H12,H14}.
For recent applications of TTs as value function approximators, see e.g.~\cite{OS21,DKK19}.
Solution methods based on high-dimensional polynomials and tensor spaces have also been considered in~\cite{KK18,DK21}. As a conjecture for future work,  block sparsity of the TTs appearing in optimal control methods could be exploited to further reduce the sample complexity \cite{trunschkeblock}.

In contrast to the aforementioned methods, our new approach is a \textit{dynamical low rank approximation} (DLRA)~\cite{KL07,KL10} of the value function.
The main idea is to approximate solutions to matrix- or tensor-valued ordinary differential equations (ODEs) by projecting the right-hand side onto the tangent space of the manifold of matrices/tensors of fixed (TT-)rank at the current approximation.
In this abstract setting, the projection is usually decomposed into orthogonal parts of the tangent space after which a splitting scheme is applied, resulting in so called projector-splitting schemes~\cite{LRS13,LOV15,KLW16,LC20,LC21}.
The obtained approximation is quasi-optimal on a finite time domain, a property known as the Dirac-Frenkel variational principle, or Dirac-Frenkel/McLachlan variational principle~\cite{F35,ML64}.
DLR approximations to parabolic PDEs have been studied in~\cite{C20,B21}, but -- to the best of our knowledge -- this work is the first application of DLR methods to a finite horizon optimal control problem and in particular to the nonlinear HJB equation.
In order to derive an abstract DLR problem on the TT manifold, we use a Variational Monte Carlo (VMC) approach ~\cite{EST20,EST21}.
In our setting it can be understood as an empirical least squares tensor regression based on random samples.

\section{The optimal control problem}
\label{sec:oc}

Throughout this work, we consider a deterministic dynamical system
\begin{align}
    \dot{x}(t) &= f(t,x(t)) + g(t,x(t))u(t), \quad t\in [t_0,T] \label{eq:lin_dyn},\\
    x(t_0) &= x_0, \label{eq:lin_ic}
\end{align}
with initial time $t_0\in [0,T]$, initial condition $x_0\in\Omega \subset \mathbb{R}^d$, control $u\in L^2(0,T;\mathbb{R}^m)$, free dynamics $f:[0,T]\times \Omega \rightarrow \Omega$ and control interface $g:[0,T]\times\Omega\rightarrow \mathbb{R}^{d\times m}$.
To ensure existence and uniqueness of solutions (for admissible controls $u$), we assume $f$ and $g$ to be smooth (possibly nonlinear) functions.
A total cost is associated with the triple $(t_0,x_0,u) \in [0,T]\times\Omega\times L^2(0,T;\mathbb{R}^m)$ in terms of the \textit{cost functional}
\begin{align}
    \mathcal{J}(t_0,x_0,u) = \int_{t_0}^T c(t ,x(t )) +  u(t)^{\intercal}R(t) u(t ) \mathrm{d}t  + c_T(x(T)), \label{eq:cost_function}
\end{align}
where the running cost $c:[0,T]\times\Omega\rightarrow \mathbb{R}_+$ and the terminal cost $c_T:\Omega\rightarrow \mathbb{R}_+$ are non-negative, coercive and smooth functionals.
Moreover, $R:[0,T]\rightarrow \mathbb{R}^{m\times m}$ is continuous, $R(t)$ is positive definite for all $t$, and the trajectory $x(\cdot)$ is subject to~\eqref{eq:lin_dyn}+\eqref{eq:lin_ic} with the given control $u$.
The function mapping time-state pairs to optimal future costs is called the \textit{value function}.
It is canonically defined as
\begin{align*}
    V^*: [0,T]\times \Omega \rightarrow \mathbb{R}, \quad (t_0,x_0) \mapsto \inf_{u\in L^2(0,T;\mathbb{R}^m)}\mathcal{J}(t_0,x_0,u).
\end{align*}
If the dynamics and cost terms satisfy sufficient regularity conditions, the value function is given as the viscosity solution of the well known Hamilton-Jacobi-Bellman equation.
\begin{theorem}[see e.g. \cite{BC97,BD97}]
\label{thm:visc_sol} 
Let $\ell(t,x,u) = c(t,x) + u^{\intercal}R(t)u$ and assume there are $\sigma, \delta \geq 1$ with $\sigma < \delta$, $\ell_0 > 0$.
Moreover, for every compact $K\subset \mathbb{R}^d$ there exists some $f_K > 0$ such that
\begin{align*}
    \| f(x,u) \| &\leq f_K(1+\| x \|^{\sigma}), \quad \textnormal{ for all }~ (x,u)\in K\times \mathbb{R}^m, \\
    |\ell(x,u)| & \geq \ell_0 \| a \|^{\delta} \quad \textnormal{ for all }~ (x,u)\in \mathbb{R}^d \times \mathbb{R}^m.
\end{align*}
Then, the value function $V^*$ is the unique viscosity solution of the HJB equation
\begin{align}
    \dfrac{\partial}{\partial t}V^*(t ,x) + \min_{u\in\mathbb{R}^m}\left[  \nabla_x V^*(t ,x)^{\intercal}(f(t ,x) + g(t ,x)u) + c(t ,x) + u^{\intercal}R(t)u \right] &= 0 \label{eq:hjb}\\
    V^*(T,\cdot) &= c_T(\cdot). \label{eq:hjb_ic}
\end{align}
\end{theorem}
Note that the HJB equation is an infinitesimal version of the Bellman equation, which we state for the sake of completeness.
\begin{theorem}[\cite{BC97,BD97}]
For all $x_0\in\Omega$ and $0\leq t_0 \leq t_1 \leq T $, we have
\begin{align}
    V^*(t_0,x_0) = \inf_{u\in L^2(t_0,t_1;\mathbb{R}^m)}\left[ \int_{t_0}^{t_1} \ell(t,x(t),u(t)) \mathrm{d}t  + V^*(t_1,x(t_1))\right], \label{eq:bellman_eq}
\end{align}
where $x(t)$ satisfies \eqref{eq:lin_dyn} with initial condition $x(t_0) = x_0$ and control $u$.
\end{theorem}
Now consider feedback controls of the form $u(t) = \alpha(t,x(t))$, where $\alpha: [0,T]\times\Omega \rightarrow \mathbb{R}^m$ is continuous on $[0,T]$ and Lipschitz in $\Omega$.
We call such functions $\alpha$ \textit{admissible feedback laws} (or equivalently \textit{admissible policies}) and we denote the set of admissible policies by $\mathcal{A}$.
Next, we define the policy evaluation function $\mathcal{J}_{\alpha}$ via the associated cost
\begin{align}
    \mathcal{J}^{\alpha}(t_0,x_0) 
    = \int_{t_0}^T c(t,x(t)) + \alpha(t,x(t))^\intercal R(t) \alpha(t,x(t)) \mathrm{d}t + c_T(x(T)). \label{eq:policy_eval}
\end{align}
An \textit{optimal policy} $\alpha^*$ is a policy which achieves minimal costs for any starting values, i.e.
\begin{align*}
    \mathcal{J}^{\alpha^*}(t_0,x_0) = \min_{\alpha\in \mathcal{A}} \mathcal{J}^{\alpha}(t_0,x_0), \quad \textnormal{ for all }~ (t_0,x_0)\in [0,T]\times \Omega.
\end{align*}
The goal of optimal (feedback) control is to approximate such an optimal policy.
If the value function is known and partially differentiable, an optimal policy can be obtained immediately, as the following theorem states.
\begin{theorem}[\cite{BC97}]
An \textit{optimal policy} is given by 
\begin{align}
    \alpha^*(t ,x) = -\dfrac{1}{2}R(t )^{-1}g(t ,x)^{\intercal} \nabla_{x} V^*(t ,x), \label{eq:opt_pol}
\end{align}
if the gradient of $V^*$ exists.
\end{theorem}
Hence, the problem of synthesizing an optimal policy is the problem of finding the value function, which involves solving the HJB equation \eqref{eq:hjb} or the Bellman equation \eqref{eq:bellman_eq}.

In the following, we always assume that the conditions of Theorem~\ref{thm:visc_sol} are satisfied so that that an optimal policy is given by the value function via~\eqref{eq:opt_pol}.
Hence, we can identify the value function with the policy evaluation function of that optimal policy denoted by $V^* = \mathcal{J}^{\alpha^*}$.
Our goal is to approximate the value function successively on small subintervals, moving backwards in time from $t=T$ to $t=0$.
This approach is based on Bellman's principle.
However, in contrast to comparable recent work~\cite{OS21}, we use the HJB equation~\eqref{eq:hjb} on each subinterval instead of the Bellman equation.
In particular, we define suitable approximate solutions to the HJB equation by means of the Dirac-Frenkel variational principle.
While theoretical simplicity is lost to some extend, computational simplicity is gained in return.
This is mainly because DLR approximations of~\eqref{eq:hjb} can be computed very efficiently since samples do not have to be propagated through the dynamics to evaluate the integral in~\eqref{eq:bellman_eq}.


Assume now that $0 = t_0 < t_1 < \hdots t_m = T$ is an equidistant discretisation and consider a partitioning $\{ [t_i,t_{i+1}] \}_{i=0}^m$ of the time interval $[0,T]$.
An immediate consequence of Bellman's principle is that an optimal policy for the whole time domain $[0,T]$ must also be optimal on any subinterval $[t_i,t_{i+1}]$.
Conversely, a policy that is optimal on all subintervals is also optimal on the whole interval.
This enables to learn the value function by moving backwards in time and (approximately) computing the restrictions $V^*(t,x)\big|_{t\in[t_i,t_{i+1}]}$ for $i=m-1,\hdots,0$.
In the following, we denote by $V^*_i = V^*\big|_{[t_i,t_{i+1}]}$ for $i=0,\hdots,m-1$ the restrictions of the value function to a particular subinterval and set $V^*_m = c_T$.
Approximations of $V^*_i$ are denoted by $\hat{V}_i$ and the approximation $\hat{V}$ of $V^*$ on the whole time domain is defined by $\hat{V}=\hat{V}_i$ on $[t_i,t_{i+1}]$.
Algorithm~\ref{alg:bellman_princ} summarizes the idea of successive backward approximation, which we deploy to approximate the value function.

\begin{algorithm}
\caption{Bellman-based backwards scheme to approximate the value function}
\label{alg:bellman_princ}
\KwData{ Time discretisation points $0=t_0 < \hdots < t_m = T$, approximation $\hat{V}_m$ of the terminal cost.}
\KwResult{Approximation $\hat{V}$ of the value function.}
\For{$i=m-1,m-2,\hdots,0$}{
    Compute approximate solution $\hat{V}_{i}$ of the HJB-eq. \eqref{eq:hjb} on $[t_i,t_{i+1}]$ with terminal condition $\hat{V}_{i+1}$. \\
     Set $\hat{V} = \hat{V}_i$ on $[t_i,t_{i+1}]$.}
\end{algorithm}
 
TT approximations of the value function by means of such a backwards scheme were already presented e.g. in~\cite{OS21}.
In that work however, the integral formulation~\eqref{eq:bellman_eq} is used exclusively, sampling trajectories $x(t)$ for given controls and adding up the costs.
In contrast, the DLR approximation method used here allows to directly work with the HJB equation~\eqref{eq:hjb}.

\section{Tensor trains as function approximators}
\label{sec:tt}

For practical computations, the approximations $\hat{V}_i$ from Algorithm~\ref{alg:bellman_princ} have to be confined to a finite-dimensional functions space.
To this end, consider a set of one-dimensional basis functions $\phi_1,\hdots,\phi_n : \mathbb{R}\rightarrow \mathbb{R}$ and functions $v:\mathbb{R}^d\rightarrow\mathbb{R}$ of the form
\begin{align}
    v(x) = A\phi(x)= \sum_{i_1,\hdots,i_d = 1}^n A_{i_1,\hdots,i_d} \phi_{i_1}(x_1)\cdot \hdots \cdot \phi_{i_d}(x_d) , \label{eq:tensor_func}
\end{align}
with coefficient tensor $A\in \mathbb{R}^{n\times n\times\hdots \times n}$ of order $d$.
Usually, the basis functions $\phi_i$ are (orthonormal) polynomials. Consequently, $v$ is a multivariate polynomial with a storage complexity of $\mathcal{O}(n^d)$ for its coefficient tensor.
The TT format provides a possibility to alleviate this exponential complexity by assuming some low-rank structure.
A TT representation of $A$ is any decomposition of the form 
\begin{align}\label{eq:tt_decomp}
    A_{i_1,\hdots,i_d} = U^{1}_{i_1}\cdot \hdots \cdot U^{d}_{i_d}, 
    \end{align}
where
\begin{align*}
    U^1\in \mathbb{R}^{n\times r_1},~ U^{\mu} \in \mathbb{R}^{r_{\mu-1}\times n \times r_{\mu}} ~\textnormal{for}~ i=2,\hdots,d-1,~ U^d\in \mathbb{R}^{r_{d-1}\times n}
\end{align*}
are called the \textit{components} of the representation and $i_{\mu}$ denotes the middle index of the component, i.e. $U^{\mu}_{i_{\mu}}\in\mathbb{R}^{r_{\mu-1}\times r_{\mu}}$.
The rank of the specific representation is given by the tuple $(r_1,\hdots,r_{d-1})$.
The TT-rank $\mathbf{r} = (r_1,\hdots,r_{d-1})$ of $A$ is defined as the (entry-wise) minimal rank tuple such that a TT representation~\eqref{eq:tt_decomp} with the corresponding ranks exists.
Such a minimal TT representation exists for any tensor.
In fact, the minimal rank entry $r_{\mu}$ is equal to the matrix rank of the $\mu$-th unfolding of $A$ (for details we refer to~\cite{HRS12}).
The TT representation exhibits a storage complexity of $\mathcal{O}(n d \max(r_1,\dots,r_{d-1})^2)$, scaling only linearly in the dimension $d$, and hence avoiding the curse of dimensionality, provided that the ranks stay bounded.
It is important to note that even for fixed rank $\mathbf{r}$, a decomposition of the form~\eqref{eq:tt_decomp} is not unique.
For any $\mu=1,\hdots,d-1$ we can set $U_{\mu} \rightarrow U_{\mu}S$ and $U_{\mu+1}\rightarrow S^{-1}U_{\mu+1}$ for invertible $S\in\mathbb{R}^{r_{\mu}\times r_{\mu}}$ without changing the tensor.
A unique representation is then given by requiring \textit{left-} and \textit{right-orthogonality} of the components in the sense of the following definition.
\begin{definition}
For a component $U_{\mu}\in\mathbb{R}^{r_{\mu-1}\times n \times r_{\mu}}$, define the \textit{left and right unfolding}
\begin{align*}
    L(U^{\mu}) \in \mathbb{R}^{r_{\mu-1}n\times r_{\mu}}, \qquad R(U^{\mu}) \in \mathbb{R}^{r_{\mu-1}\times r_{\mu}n},
\end{align*}
by suitable matrix reshaping (for details regarding the order, see e.g.~\cite{Steinlechner16}).
A component $U^{\mu}$ is called left- or right-orthogonal if
\begin{align*}
    L(U^{\mu})^{\intercal}L(U^{\mu}) = I_d \in \mathbb{R}^{r_{\mu}\times r_{\mu}}, \quad \textnormal{or} \quad R(U^{\mu})R(U^{\mu})^{\intercal} = I_d \in \mathbb{R}^{r_{\mu-1}\times r_{\mu-1}},
\end{align*}
respectively.
A TT representation $A_{i_1,\hdots,i_d} = U^{1}_{i_1}\cdot \hdots \cdot U^{d}_{i_d}$ of a tensor $A$ is called $\mu$-orthogonal if $U^1,\hdots, U^{\mu-1}$ are left orthogonal and $U^{\mu+1},\hdots,U^{d}$ are right orthogonal.
In that case, $U^{\mu}$ is called the \textit{core of the representation}.
\end{definition}

Left and right orthogonality of all but one component imposes $\sum_{\mu=1}^{d-1} r_{\mu}^2$ additional conditions on the representation.
Hence, the $\mu$-orthogonal TT representation of $A$ is unique for any $\mu$.

For a given TT rank $\mathbf{r}$, we define the set 
\begin{align*}
    \mathcal{M}_{\mathbf{r}} = \{ A\in\mathbb{R}^{n\times\hdots\times n} : A~\textnormal{ has TT rank }\mathbf{r} \}.
\end{align*}
 It is noteworthy that $\mathcal{M}_{\mathbf{r}}$ is a smooth manifold in $\mathbb{R}^{n\times \hdots \times n}$~\cite{HRS12}.
 With a chosen suitable basis $\{ \phi_1,\hdots, \phi_n \}$, we define a set of function approximations
 \begin{align*}
     \mathfrak{F}_{\mathbf{r}} = \left\{ v:\mathbb{R}^d\rightarrow \mathbb{R}: v ~\textnormal{admits a representation}~ \eqref{eq:tensor_func} ~\textnormal{with}~ A\in \mathcal{M}_{\mathbf{r}} \right\}.
 \end{align*}
Note that by identification of a function with its coefficient tensor, $\mathfrak{F}_{\mathbf{r}}$ forms a smooth manifold in the $n^d$-dimensional linear space $\mathfrak{F} = \bigcup_{\mathbf{r}}\mathfrak{F}_{\mathbf{r}}$ in the same way that $\mathcal{M}_{\mathbf{r}}$ forms a smooth manifold in $\mathbb{R}^{n\times\hdots\times n}$.
In order to do perform an optimisation on $\mathfrak{F}_{\mathbf{r}}$, or $\mathcal{M}_{\mathbf{r}}$, respectively, we require a representation of the tangent space $\mathcal{T}_U(\mathcal{M}_{\mathbf{r}})$ of $\mathcal{M}_{\mathbf{r}}$ in $U$.
Throughout this work, we use the following representation.

\begin{theorem}[\cite{HRS12} or \cite{Steinlechner16}]
\label{thm:tangent_space}
Let $U\in \mathcal{M}_{\mathbf{r}}$ be $d$-orthogonal.
The \textit{tangent space} $\mathcal{T}_{U}(\mathcal{M}_{\mathbf{r}})$ of $\mathcal{M}_{\mathbf{r}}$ in the point $U$ is given by $\tau(X)$, where 
\begin{align*}
    X &= U^\ell_1 \times \hdots \times U^\ell_{d-1} \times \mathbb{R}^{r_{d-1}\times n\times r_d}, \\
    U_{\mu}^\ell &= \{ W^{\mu}\in \mathbb{R}^{r_{\mu-1}\times n\times r_{\mu}} : L(U^{\mu})^\intercal L(W^{\mu}) = 0 \in \mathbb{R}^{r_{\mu}\times r_{\mu}} \},
\end{align*}
and
\begin{align}
    &\tau: X \longrightarrow \mathcal{T}_U(\mathcal{M}_{\mathbf{r}}), \quad \tau(W^1,\hdots,W^d) = \delta U \notag\\
    &\delta U_{i_1,\hdots,i_d} = \sum_{\mu=1}^d U^1_{i_1}\cdot\hdots\cdot U^{\mu-1}_{i_{\mu-1}}W^{\mu}_{i_{\mu}} U^{\mu+1}_{i_{\mu+1}}\cdot\hdots \cdot U^d_{i_d}. \label{eq:tangent_rep}
\end{align}
\end{theorem}
The tangent space has the same dimension as the underlying manifold.
The previously mentioned ambiguity in the representation is now eliminated due to the \textit{gauging conditions} $L(U^{\mu})^{\intercal} L(W^{\mu})=0$ in $U^{\ell}_{\mu}$.
\begin{corollary}
Each of the spaces $U_{\mu}^\ell$ has dimension $r_{\mu-1}n r_{\mu} - r_{\mu}^2$ and hence the tangent space has dimension
\begin{align*}
    n_X \coloneqq \textnormal{dim}(X) = \sum_{\mu=1}^d r_{\mu-1}n r_{\mu} - \sum_{\mu=1}^{d-1} r_{\mu}^2.
\end{align*}
\end{corollary}


Using the representation~\eqref{eq:tangent_rep} for elements $\delta U$ of the tangent space of $U$, a simple form for the sum $U + \delta U$ can be obtained.

\begin{lemma}[see~\cite{Steinlechner16}]
\label{lem:tangent_add}
Let $U\in \mathcal{M}_{\mathbf{r}}$ be $d$-orthogonal and denote its component tensors by $U_1,\hdots,U_d$.
Let $\delta U \in \mathcal{T}_{U}(\mathcal{M}_{\mathbf{r}})$ be given by $(\delta U_1,\hdots, \delta U_d) \in X$.
Then, 
\begin{align*}
U(i_1,\hdots,i_d) + \delta U(i_1,\hdots,i_d) = 
     &\begin{bmatrix}  \delta U_1(i_1) & U_1(i_1) \end{bmatrix} \begin{bmatrix} U_{2}(i_{2}) & 0 \\ \delta U_{2} (i_{2}) & U_{2}(i_{2}) \end{bmatrix} \hdots \\
     &\hdots \begin{bmatrix} U_{d-1}(i_{d-1}) & 0 \\  \delta U_{d-1} (i_{d-1}) & U_{d-1}(i_{d-1}) \end{bmatrix}
     \begin{bmatrix} U_d(i_d) \\ U_d(i_d) + \delta U_d(i_d) \end{bmatrix}.
\end{align*}
\end{lemma}
This can easily be verified by multiplying out the matrix products.
In particular, the sum $U+\delta U $ has at most TT-rank $2\mathbf{r}$. 

\section{The Dirac-Frenkel variational principle}
\label{sec:dlra}

The Dirac-Frenkel variational principle~\cite{F35} provides a principled way to approximate tensor valued ODEs of the form
\begin{align}
    \dot{A}(t) &= F(t,A(t)),\label{eq:tensor_ode}\\
    A(0) &= A_0, \label{eq:tensor_ode_ic}
\end{align}
where $A(t)\in\mathbb{R}^{n\times\hdots\times n}$, on the manifold $\mathcal{M}_{\mathbf{r}}$.
More precisely, given an approximation $Y_0 \in \mathcal{M}_{\mathbf{r}}$ of the initial condition $A_0$, an approximation $Y(t)\in\mathcal{M}_{\mathbf{r}}$ of $A(t)$ is defined as the solution of the TT-valued ODE
\begin{align}
    \dot{Y}(t) &= \argmin_{\vartheta\in\mathcal{T}_{\mathcal{M}_\mathbf{r}}(Y(t))}\| \vartheta - F(t,Y(t)) \|, \label{eq:tt_ode}\\
    Y(0) &= Y_0. \label{eq:tt_ode_ic}
\end{align}
The minimum in~\eqref{eq:tt_ode} is attained by the orthogonal projection of the right-hand side onto the tangent space, leading to
\begin{align}
    \dot{Y}(t) = P_{\mathcal{T}_{\mathcal{M}_r}(Y(t))} F(t,Y(t)). \label{eq:dirac_frenkel}
\end{align}
In this abstract setting, error bounds can be derived, which we quote for the sake of completeness.

\begin{theorem}\cite{LRS13}
Suppose that $\dot{A}(t) \leq \mu$ and that a continuously differentiable best approximation $X(t) \in \mathcal{M}_{\mathbf{r}}$ to $A(t)$ exists for $t\in[0,T]$.
Let $\delta > 0$ be such that the smallest nonzero singular value of every matrix unfolding of $X(t)$ is greater or equal to $\rho$, and assume that the best-approximation error is bounded by $\| X(t) - A(t) \| \leq c\rho$ for $t\in[0,T]$ with a constant $c$ depending only on the dimension $d$.
Then, the approximation error of the dynamical low-rank approximation defined by~\eqref{eq:tt_ode_abstract} with initial value $Y(0) = X(0)$ is bounded by
\begin{align*}
   \| Y(t) - X(t)\| \leq 2\beta e^{\beta t} \int_0^t
 \| X(s) - A(s)\| \mathrm{d}s,
\end{align*}
with $\beta = C\mu\rho-1$ for $t\in[0,T]$, as long as the right-hand side remains bounded by $c\rho$.
The constant $C$ is only dependent on $d$ and is given in~\cite{LRS13}.
\end{theorem}

In recent years there have been numerous works on the numerical treatment of ODEs of this type, see~\cite{KL07,LO13,KLW16} for an introduction in the matrix case and~\cite{LOV15,LC20,LC21} for more recent tensor-based research directions.
Generally, these methods rely on a splitting of the projector $P_{\mathcal{M}_{\mathbf{r}}(Y(t))}$ into orthogonal parts of the tangent space, so-called \textit{projector splitting algorithms}.
The norm $\|.\|$ governing~\eqref{eq:tt_ode} and hence the projector is usually the Frobenius norm.
This is in contrast to our work, where $\|.\|$ is an empirical norm\footnote{the details of which are provided in the next chapter}.
Carrying over results from the treatment of the abstract Dirac-Frenkel principle to the empirical case (specifically the projector splitting schemes) is an important direction of future work, that we do not yet address in this paper.

\section{Dynamical low-rank approximation of the HJB equation}
\label{sec:pi}

Based on the preceding review of tools that we require, we now return to the HJB equation~\eqref{eq:hjb} on $[t_{i},t_{i+1}]$ with terminal condition $\hat{V}_{i+1}(t_{i+1},\cdot)$.
The goal is to obtain an approximation $\hat{V}_i$ of the value function on the current interval.
Inserting~\eqref{eq:opt_pol} into the HJB~\eqref{eq:hjb} leads to a coupled problem:\\[1ex]
Find $V$ such that
\begin{align}
    \dfrac{\partial}{\partial t}V(t ,x) +  \nabla_x V(t ,x)^{\intercal}(f(t ,x) + g(t ,x)\alpha(t,x)) + c(t ,x) + \alpha(t ,x)^{\intercal}R(t)\alpha(t ,x) &= 0, \label{eq:ghjb} \\
    V(t_{i+1},\cdot) &= \hat{V}_{i+1}(t_{i+1},\cdot),\label{eq:ghjb_ic} 
\end{align}
where $\alpha$ satisfies 
\begin{align}
    \alpha(t ,x) = -\dfrac{1}{2}R(t )^{-1}g(t ,x)^{\intercal} \nabla_{x} V(t ,x). \label{eq:pol_upd}
\end{align} 
To compute $\hat{V}_{i}$, we use a fixed point iteration of the coupled problem, iteratively solving~\eqref{eq:ghjb}+\eqref{eq:ghjb_ic} for fixed $\alpha$ and then updating $\alpha$ via~\eqref{eq:pol_upd}.
This procedure is known as \textit{policy iteration} in the optimal control literature.
We depict a conceptual summary in algorithm~\ref{alg:pi}.
If the solutions to~\eqref{eq:ghjb} are exact, it converges under mild assumptions on dynamics and cost terms~\cite{SL79}.
In order to track the convergence of the scheme under approximations, we introduce on $L^2((t_i,t_{i+1});L^2(\Omega,\rho))$ and $L^2((0,T);L^2(\Omega,\rho))$ the norms
\begin{align*}
    \| v \|_i^2 = \int_{t_i}^{t_{i+1}} \|v(t ,\cdot)\|^2_{L^2(\Omega,\rho)} \mathrm{d}t ,
    \qquad  \| v \|^2 = \sum_{i=0}^{m-1} \|v\|^2_i,
\end{align*}
and stop the iteration once the $\| . \|_i$-difference of two consecutive approximations becomes smaller than a specified threshold.

\begin{algorithm}
\caption{Policy iteration on subinterval}\label{alg:pi}
\KwData{ Interval $[t_i,t_{i+1}]$, terminal condition $\hat{V}_{i+1}$, admissible policy $\alpha$, error tolerance $\delta$.}
\KwResult{ Approximation $\hat{V}_i$. } 
\While{\textnormal{norm change of} $\hat{V}_i > \delta$}{
    Compute approximate solution $\hat{V}_i$ of \eqref{eq:ghjb}+\eqref{eq:ghjb_ic} on $[t_i,t_{i+1}]$. \\
    Update $\alpha \propto \nabla_x \hat{V}_i$ according to \eqref{eq:pol_upd}.
    }
\end{algorithm}

It remains to be shown how to compute the approximations $\hat{V}_i$.
To ease notation, without loss of generality we consider the interval $[t_0,t_1]$ instead of $[t_i,t_{i+1}]$ for the remainder of this chapter.
We construct $\hat{V}_0$ as a dynamical low-rank approximation of~\eqref{eq:ghjb} in the tensor train format.

Let the terminal condition $\hat{V}_{1}\in \mathfrak{F}_{\mathbf{r}}$ and consider for given $\alpha$ the following problem:\\[1ex]
Find $V$ such that
\begin{align}
    \dfrac{\partial}{\partial t}V(t,\cdot) &= \argmin_{\vartheta\in\mathcal{T}_{V(t,\cdot)}(\mathfrak{F}_{\mathbf{r}})} \| \vartheta - \tilde{F}(t,V(t,\cdot)) \|^2_{L^2(\Omega,\rho)}, \label{eq:tt_ode_abstract}\\
     V(t_{1},\cdot) &= \hat{V}_{1}(t_{1},\cdot), \label{eq:tt_ode_ic_abstract}
\end{align}
where $\tilde{F}(t,V(t,\cdot))(\cdot) = -\nabla_x V(t,\cdot)^{\intercal}(f(t,\cdot) + g(t,\cdot)\alpha(t,\cdot)) - c(t,\cdot) + \alpha(t,\cdot)^{\intercal}R(t)\alpha(t\cdot)$.
Note that this essentially means that the time derivative of $V_0$ is approximated in the tangent space of the current solution.
By a simple time inversion $t\rightarrow t_0 + (t_1-t)$, the terminal condition can be turned into an initial condition.
Crucially, any solution to~\eqref{eq:tt_ode_abstract} stays on the manifold $\mathfrak{F}_{\mathbf{r}}$ and can therefore be identified with a time-dependent coefficient tensor $A(t)\in\mathcal{M}_{\mathbf{r}}$ via $V(t,x) = A(t)\phi(x)$.
Denoting the coefficient tensor of $\hat{V}_{1}$ by $\hat{A}_1$, we see that the abstract problem~\eqref{eq:tt_ode_abstract}+\eqref{eq:tt_ode_ic_abstract} is equivalent to the TT-valued ODE
\begin{align}
    \dot{A}(t) &= \argmin_{B\in\mathcal{T}_{A(t)}(\mathcal{M}_{\mathbf{r}})} \| B\phi - F(t,A(t)\phi) \|^2_{L^2(\Omega,\rho)}, \label{eq:tt_ode_explicit}\\
    A(t_0) &= A_1, \label{eq:tt_ode_ic_explicit}
\end{align}
where $F$ arises from time inversion of $\tilde{F}$.

In general, the $L^2$-integral on the right-hand side of~\eqref{eq:tt_ode_explicit} is difficult to compute.
Nevertheless, we can easily carry out a pointwise evaluation of the basis functions $\phi$ as well as the other terms in $F(t,A(t)\phi)$.
In practice, we hence replace the exact $L^2$-norm with a Monte Carlo approximation
\begin{align*}
    \| v \|^2_{L^2(\Omega,\rho,M)} = \dfrac{1}{M}\sum_{k=1}^M |v(x_k)|^2, \qquad x_k \sim \rho,
\end{align*}
for $v\in L^2(\Omega,\rho)$.
This turns the right-hand side of the ODE into an empirical risk minimisation.
We eventually arrive at 
\begin{align}
    \dot{A}(t) &= \argmin_{B\in\mathcal{T}_{A(t)}(\mathcal{M}_{\mathbf{r}})}\dfrac{1}{M}\sum_{k=1}^M | B\phi(x_k) - F(t,A(t)\phi(x_k))(x_k) |^2, \qquad t\in(t_0,t_1)\label{eq:emp_risk_min}\\
    A(t_0) &= \hat{A}_1. \label{eq:emp_risk_min_ic}
\end{align}
Statistical bounds for the error of the empirical minimiser in~\eqref{eq:emp_risk_min} compared to the best $L^2$-approximation $\Phi^*(t) = \argmin_{\Phi\in L^2(\Omega,\rho)} \| \Phi -F(t,A(t)\phi(\cdot))(\cdot)\|_{L^2(\Omega,\rho)}^2$ are given in~\cite{EST19}.

A crucial observation is that the minimisation on the right-hand side is a linear problem since the optimisation is over the linear tangent space.
Implementation details on how the minimum in~\eqref{eq:emp_risk_min} for a given $t$ can be computed are given in Appendix~\ref{sec:tangent_fit}.
Since the fit is linear, issues of local minima are avoided which for instance occur in the alternating linear scheme (ALS)~\cite{HRS12} and other nonlinear optimisation methods.
Alternating methods can still be applied here to divide the problem into smaller sub-problems and reduce the computational burden, leading (in their simplest form) to an effective Lie-Trotter type splitting of the right hand side.
A more detailed examination of this topic is however beyond the scope of this paper and might be addressed in future work.

The numerical realisation of \eqref{eq:emp_risk_min} poses an additional hurdle.
While the true solution $A(t)$ always stays on the manifold $\mathcal{M}_{\mathbf{r}}$, it is straightforward to see that any one step with a numerical integrator, e.g. a Runge-Kutta method, leads to leaving it.
This is due to the fact that by Lemma~\ref{lem:tangent_add} any sum $U+\delta U$ where $U\in\mathcal{M}_{\mathbf{r}}$ and $\delta U \in \mathcal{T}_{U}(\mathcal{M}_{\mathbf{r}})$ has rank $2\mathbf{r}$ in general.
We therefore need to retract back onto the manifold after each step of the integrator by truncating the ranks appropriately.
To make this precise, let $t_0=t^{(0)} < t^{(1)} < t^{(2)} < \hdots < t^{(L)} = t_1$ be a \textit{micro-discretisation} of the \textit{macro-interval} $[t_0,t_1]$ with equidistant step size $\tau$ and define a numerical approximation $A_\ell$ of $A(t^{(\ell)})$ by the explicit Euler scheme
\begin{align*}
    A_0 &= \hat{A}_1, \\
    A_{\ell+1} &= \mathcal{R}(A_\ell + \tau \Delta A_\ell), \qquad \ell=0\hdots,L-1.
\end{align*}
Here, $\Delta A_\ell$ is the solution to the minimisation problem on the right-hand side of~\eqref{eq:emp_risk_min} if $A_\ell$ is substituted for $A(t)$, the addition $A_\ell+\tau\Delta A_\ell$ is performed like in Lemma~\ref{lem:tangent_add}, and $\mathcal{R}$ denotes the rank-truncation of a TT with rank $2\mathbf{r}$ back to a tensor of rank $\mathbf{r}$.
This truncation is performed by a TT-SVD with fixed rank~\cite{O09}.
Once all $A_\ell$ are obtained in this way, we define $\hat{V}_0(t,x)$ by linear interpolation, i.e.
\begin{align*}
    \hat{V}_0(t,x) = A_\ell\phi(x) + \dfrac{t-t^{(\ell)}}{\tau} (A_{\ell+1}-A_\ell)\phi(x) \quad \textnormal{for} \quad t \in [t^{(\ell)},t^{(\ell+1)}],
\end{align*}
or by simply always setting it to
\begin{align*}
    \hat{V}_0(t,x) = A_\ell\phi(x) \quad \textnormal{for} \quad t \in [t^{(\ell)},t^{(\ell+1)}].
\end{align*}

Now, if $\hat{V}^{\textnormal{old}}_0$ is the approximation from the previous policy iteration step, one could compute the empirical approximation to the $\|.\|_0$-norm
\begin{align*}
    \| \hat{V}_{0} - \hat{V}^{\textnormal{old}}_0 \|_{0,L,M}^2 &= \dfrac{1}{L-1}\sum_{\ell=0}^{L-1} \|\hat{V}_0(t^{(\ell)} ,\cdot)-\hat{V}^{\textnormal{old}}_0(t^{(\ell)},\cdot)\|^2_{L^2(\Omega,\rho,M)} \\
    &= \dfrac{1}{(L-1)M}\sum_{\ell=0}^{L-1}\sum_{k=1}^M |\hat{V}_0(t^{(\ell)} ,x_k)-\hat{V}^{\textnormal{old}}_0(t^{(\ell)},x_k)|^2
\end{align*}
and stop the iteration once this norm difference becomes smaller than the threshold $\delta$.
However, since we are first and foremost interested in obtaining a nearly optimal control $\alpha$, we instead add the change in the controls $\alpha\propto \nabla_x \hat{V}_0$ and $\alpha^{\textnormal{old}}\propto \nabla_x\hat{V}^{\textnormal{old}}_0$ and stop the iteration once
\begin{align}\label{eq:pi_conv_crit}
    \dfrac{1}{(L-1)M}\sum_{\ell=0}^{L-1}\sum_{k=1}^M |\hat{V}_0(t^{(\ell)} ,x_k)-\hat{V}^{\textnormal{old}}_0(t^{(\ell)},x_k)|^2 + \|\alpha(t^{(\ell)} ,x_k)-\alpha^{\textnormal{old}}(t^{(\ell)},x_k)\|^2 < \delta.
\end{align}
The reason for this is that the $L^2$-norm is agnostic to errors in the gradients, which may arise due to overfitting.
By requiring~\eqref{eq:pi_conv_crit}, we demand that not only $\hat{V}_0$ but also the relevant part of the gradient $\nabla_x\hat{V}_0$ converges.
In that sense, the left-hand side of~\eqref{eq:pi_conv_crit} can be seen as an empirical approximation of an $H^1$-norm of $\hat{V}_0-\hat{V}_0^{\textnormal{old}}$, where the norms for the gradients are now weighted by $R$ and $g$ to represent only the gradient parts relevant for the control.

\section{Numerical tests}
\label{sec:num}

This chapter is concerned with numerical experiments that illustrate the performance of the proposed DLR approximation\footnote{All computations are carried out on an Intel Xeon Gold 6154 CPU 3.00GHz, openSUSE Leap 15.2 distribution.}.
We consider a problem of the form
\begin{align*}
    \dot{x} = Ax + \mathrm{nl}(x) +gu,
\end{align*}
where $x\in\mathbb{R}^d$, $A\in\mathbb{R}^{d\times d}$, $g\in \mathbb{R}^{d}$, $u$ is scalar and $\mathrm{nl}$ is a smooth nonlinear function with $\mathrm{nl}(0)=0$.
In particular, the optimal control problem is derived from a modified one dimensional heat equation
\begin{alignat*}{3}
&\dfrac{\partial}{\partial t}{x}(s,t) &&= \sigma\dfrac{\partial^2}{\partial s^2} x(s,t) + x(s,t)^3 + g(s)u(t), \quad &&\text{ for } (s,t)\in[-1,1]\times (0,T), \\
& x(s,0) &&= \tilde{x}_0(s), \quad &&\text{ for } s \in [-1,1],\\
& \dfrac{\partial}{\partial s}x(-1,t) &&= \dfrac{\partial}{\partial s}x(1,t) = 0 \quad &&\text{ for } t \in (0,T),
\end{alignat*}
with unstable reaction term $x(s,t)^3$, diffusion $\sigma > 0$, scalar control $u$ and initial state $\tilde{x}_0$.
Note that due to the instability introduced by the reaction term, this problem is generally more difficult to control than most other canonically treated examples like viscous Burgers' type equations, Allen-Kahn or degenerate Zeldovich equations~\cite{KK18,OS21} since the quadratic regulator usually provides a strong and mostly stable controller for these types of problems.
This however is not the case for the nonlinear reaction problem defined above.
We hence omit the mentioned alternative examples and just note that our method can be applied with them as well, although the difference to the linear quadratic regulator would turn out to be small.

Our goal is to find a control $u$ such that the quadratic cost functional
\begin{align*}
\tilde{J}(0,x_0,u)= \int_0^{T} \| x(\cdot,t) \|_{L^2(\Omega)}^2 + \gamma u(t)^2 ~ \mathrm{d}t + c_T\|x(\cdot,T)\|_{L^2(\Omega)},
\end{align*}
is minimal with $\gamma, c_T > 0$.
A semi-discretisation of the PDE with finite differences at $d$ equidistant points $-1 =  s_1 < \hdots < s_{d} = 1$ leads to a an ODE of the form
\begin{alignat}{2}
&\dot{x} &&= Ax + x^3 + gu, \label{eq_probLin}\\
&x(0) &&= x_0, \label{eq_probbLinIni}
\end{alignat}
with $x_0 = (\tilde{x}_0(s_1),\hdots,\tilde{x}_0(s_{n}))^\intercal$, $x(t)\in\mathbb{R}^{d}$, $g = (g(s_1),\hdots, g(s_d))^\intercal\in\mathbb{R}^{d}$ and $A\in\mathbb{R}^{d\times d}$ is given by
\begin{align*}
A = \dfrac{\sigma}{h^2}\begin{pmatrix}
-2 & 2 & &\\
1 & -2 & 1 &\\
& \ddots & \ddots &\ddots  &\\
& & 1 & \ddots & 1 \\
& & & 2 & -2
\end{pmatrix}, \quad h = s_1-s_0 = \dfrac{2}{d-1}.
\end{align*}
The $x$-dependent term in the cost functional can be approximated using a simple quadrature rule with nodes $s_1,\hdots, s_d$ (here, we use the rectangle rule with an additional node at the last grid point $s_d$).
This yields the new cost functional
\begin{align}
J(0,x_0,u) = \int_0^{T} x(t)^\intercal Q x(t) + \gamma u(t)^2 ~ \mathrm{d}t + c_T x(T)^\intercal Q x(T),\label{probLinCost}
\end{align}
where
\begin{align*}
 Q = h\begin{pmatrix}
1 & &   \\
 & \ddots &  \\
& &  1
\end{pmatrix},
\end{align*}
and $x(t)$ is understood to be the solution of $\dot{x} = Ax + x^3 + gu$ with starting value $x_0$.
The control problem is now to find a control $u$ for the nonlinear system~\eqref{eq_probLin} such that~\eqref{probLinCost} is minimal for every starting value $x_0$. 
 
To specify the control problem, we choose the parameters $\sigma=1$, $\gamma=0.1$, $c_T = 1$ and $g = \chi_{[-0.4,0.4]}$ and discretise with $n=12$ equidistant grid points.
The time horizon is $T=0.3$ and the time step size $\tau = t_{i+1}-t_i$ is $0.001$, which is used for both the macro-intervals as well as the micro-intervals of the policy iteration (see Section~\ref{sec:pi}).
The same step size is also used to discretise the integral in~\eqref{probLinCost} when computing the costs.
As a threshold for the policy iteration, $\delta = 10^{-6}$ is set.
We choose $\Omega = (-2,2)^d$ and let $\rho$ be the uniform distribution on $\Omega$.
For the TT approximations, we use the first $n$ $H^2_{\textnormal{mix}}(\Omega)$-orthonormal polynomials $\phi_1,\hdots, \phi_n$ as basis functions (up to degree $n-1$).
Here, $H^2_{\textnormal{mix}}(\Omega)$ denotes the tensorised space $\bigotimes_{\mu=1}^d H^2((-2,2))$, $H^2$ is the Sobolev space of twice weakly differentiable functions.
We set $n=9$, yielding a maximal polynomial degree of $8$ in the basis.
The rank of the TT manifold is chosen to be
 \begin{align*}
     \mathbf{r} = (3,5,5,5,5,5,5,5,5,5,3).
 \end{align*}
Note that by this the dimension of the approximation space is reduced from $n^d = 9^{12} > 282$ trillion to a manageable number of degrees of freedom $\leq n  d  \max_{\mu=1,\hdots,d}(r_{\mu})^2 = 2700$.
The number of sample points used to approximate the $L^2$-norm in~\eqref{eq:emp_risk_min} is chosen as 
\begin{align*}
    M = 6\cdot n d  \max_{\mu=1,\hdots,d}(r_{\mu})^2 = 6\cdot 8 \cdot 12 \cdot 5^2 = 16200,
\end{align*}
which is a generous upper bound for the number of degrees of freedom of the fit.

As a benchmark for assessing the performance of our method, we use the TT-based approach from~\cite{OS21} with the same hyper-parameters.
To make this precise, instead of solving~\eqref{eq:emp_risk_min} by means of our dynamical low-rank scheme, $\hat{V}_i$ is approximated in each policy iteration step by sampling the trajectories $x_k(t)$, $t\in[t_i,t_{i+1}]$ of all sample points.
With this, the integrals
\begin{align*}
    \hat{V}_i(t_i,x_k) = \int_{t_i}^{t_{i+1}} \ell(t,x_k(t),\alpha(t,x_k(t)) \mathrm{d}t  + \hat{V}_{i+1}(t_{{i+1}},x_k(t_{i+1}))
\end{align*}
are evaluated subsequently.
An approximation of $V_i(t_i,\cdot)$ is then obtained via a nonlinear fit of a rank-$\mathbf{r}$ TT to the resulting data-target pairs $(x_k,y_k=\hat{V}_i(t_i,x_k))_{k=1}^M$, which is performed by the ALS.
Note that the authors in~\cite{OS21} suggest replacing the upper integral bound $t_{i+1}$ with $t_{i+l}$, $l>1$, where the trajectory on $[t_{i+1},t_{i+l}]$ is controlled by the already computed (nearly optimal) controls from previous steps, to remove the error associated with $\hat{V}_{i+1}$ from the computation of $\hat{V}_i$.
Since this greatly increases the computational complexity, we stick with the above mentioned \textquotedblleft one-step scheme\textquotedblright ~and refer to this benchmark method as the \textit{Bellman} method, since it explicitly utilises Bellman's principle in the form of the terminal cost $\hat{V}_{i+1}$. Our method, utilising Dynamical Low Rank Approximation, will be called the DLRA method.
Even for the DLRA method we have found it beneficial for stable convergence to compute some $\hat{V}_i$ with the Bellman method before starting the dynamical low rank solver.
In this example the first $10$ of the $300$ approximations are computed in this way.
 
\begin{remark}
 In both the nonlinear fit required for the Bellman method and the linear fit of our DLRA method, we add a regularisation term $\delta\| \hat{V}_i(t_i,\cdot) \|^2_{H^2_{\textnormal{mix}}(\Omega)}$ to the minimisation functional.
 Due to the multilinear structure of the TT and our choice of the basis functions as $H^2_{\textnormal{mix}}(\Omega)$-orthonormal, this leads to local minimisation problems of the form
 \begin{align*}
     \min_{ \mathbf{c}} \| M\mathbf{c} - \mathbf{y} \|_2^2 + \delta \| \mathbf{c} \|^2_F
 \end{align*}
 in ALS (compare to~\cite{OS19}).
 Here, $\mathbf{c}\in\mathbb{R}^{r_{\mu-1}\times n_{\mu}\times r_{\mu}}$ denotes the core that is currently optimised and $\|.\|_F$ denotes the Frobenius norm.
 In both methods, we use $\delta = 10^{-10}$ but in ALS we successively lower $\delta$ via
 \begin{align*}
     \delta \rightarrow \max(0.9,\| M\mathbf{c} - \mathbf{y} \|_2^2/ \| \mathbf{y} \|^2_2)\cdot \delta
 \end{align*}
 after every sweep.
 This is a purely heuristical rule to make sure the regularisation is relaxed once the attractor of the global minimum is found.
\end{remark}
 
 \begin{remark}
 For the DLRA method we add an additional regularisation term 
 \begin{align*}
     \delta_0 | B\phi(0) |^2
 \end{align*}
 to the minimisation in~\eqref{eq:emp_risk_min} since we know that the right-hand side satisfies $F(t,A(t)\phi(0))=0$.
 Note that this can be realised by simply adding the point $x_{M+1} = 0$ to the set of samples $\{x_k\}_{k=1}^M$.
 Since this is a hard constraint on the true solution, we set $\delta_0 = 10^{10}$.
 \end{remark}
 
As a second, classical benchmark, we consider the linear quadratic regulator (LQR), resulting from linearising the problem around $x=0$.
Since this controller does not see the unstable reaction term, we expect poor performance compared to both the Bellman and the DLRA method.
 
To compare the practical performance of the methods, two different sets of initial conditions $\tilde{x}_0$ are generated.
For the first set, we sample a polynomial degree between $2$ and $20$ and then again randomly sample the coefficients of a univariate polynomial of that degree.
Denoting this polynomial $p$, we then set $\tilde{x}_0(s) = (s-1)^2(s+1)^2p(s)$ to make sure $\tilde{x}_0$ satisfies the Neumann boundary conditions.
Finally, in order to have interesting trajectories \eqref{eq_probLin} for which the $x^3$-term requires strong control beyond LQR, we normalise such that $\max_{s\in[-1,1]} |\tilde{x}_0(s)| = 1.9$.
The second set of initial conditions is generated by simply setting $\tilde{x}_0(s) \equiv c$ for constants $c\in [1,2)$.
 
\begin{figure}
    \centering
    \includegraphics[scale=0.6]{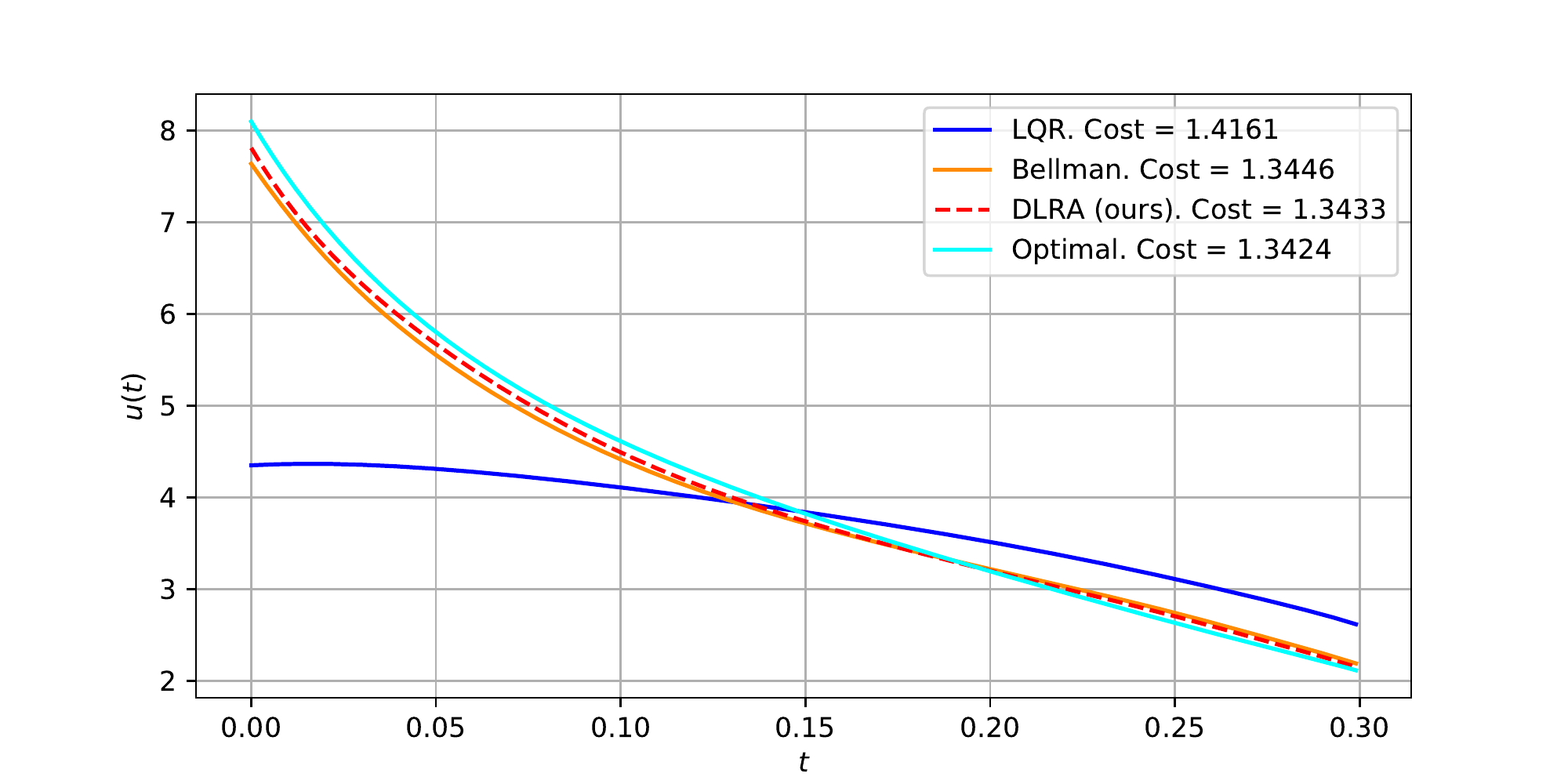}
    \caption{Control values $u(t)$ of the different controllers along the trajectory of a fixed polynomial initial condition $x_0$.}
    \label{fig:control_poly_ic}
\end{figure}

\begin{figure}
    \centering
    \includegraphics[scale=0.6]{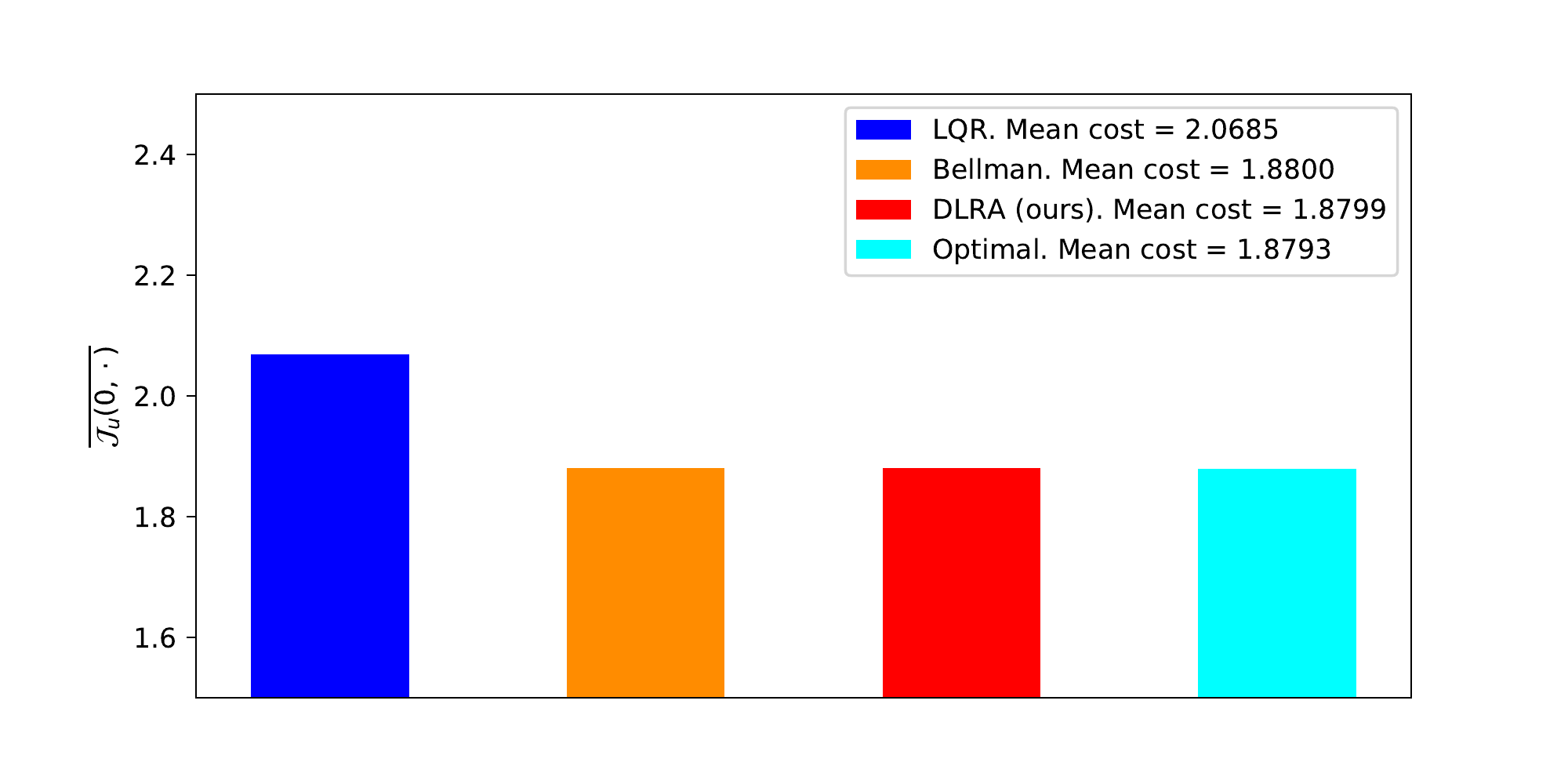}
    \caption{Sample-mean costs $\overline{\mathcal{J}_u(0,\cdot)} \approx \dfrac{1}{N}\sum_{k}^N J_u(0,x^{(k)}_0)$ of polynomial initial conditions $x_0^{(k)}$ with the different controllers.}
    \label{fig:costs_poly_ic}
\end{figure}

\begin{figure}
    \centering
    \includegraphics[scale=0.6]{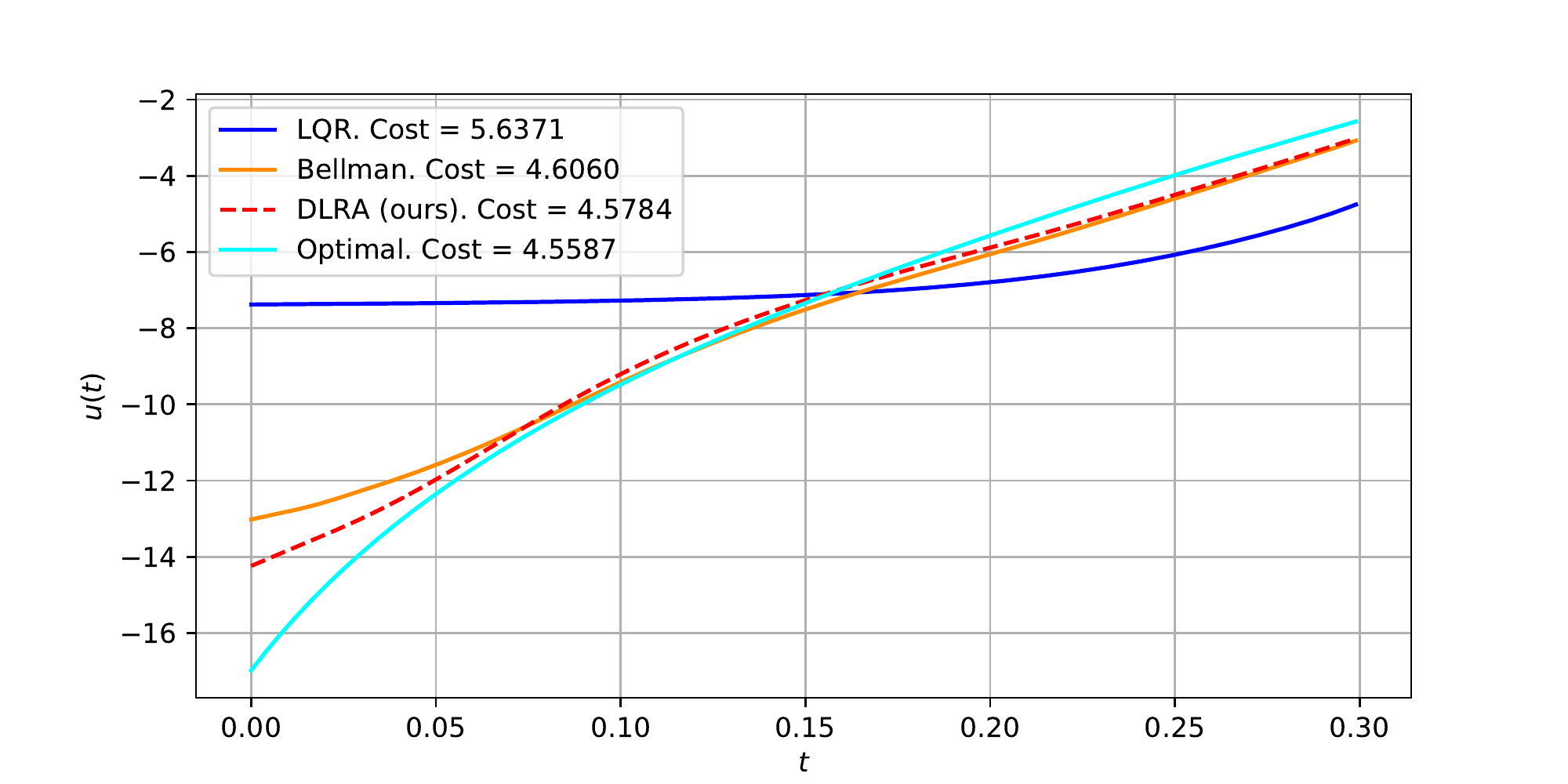}
    \caption{Control values $u(t)$ of the different controllers along the trajectory of a fixed uniform initial condition $x_0 = c\cdot (1,\hdots,1)^{\intercal}$. In this example $c=1.28$ is used.}
    \label{fig:control_const_ic}
\end{figure}

\begin{figure}
    \centering
    \includegraphics[scale=0.6]{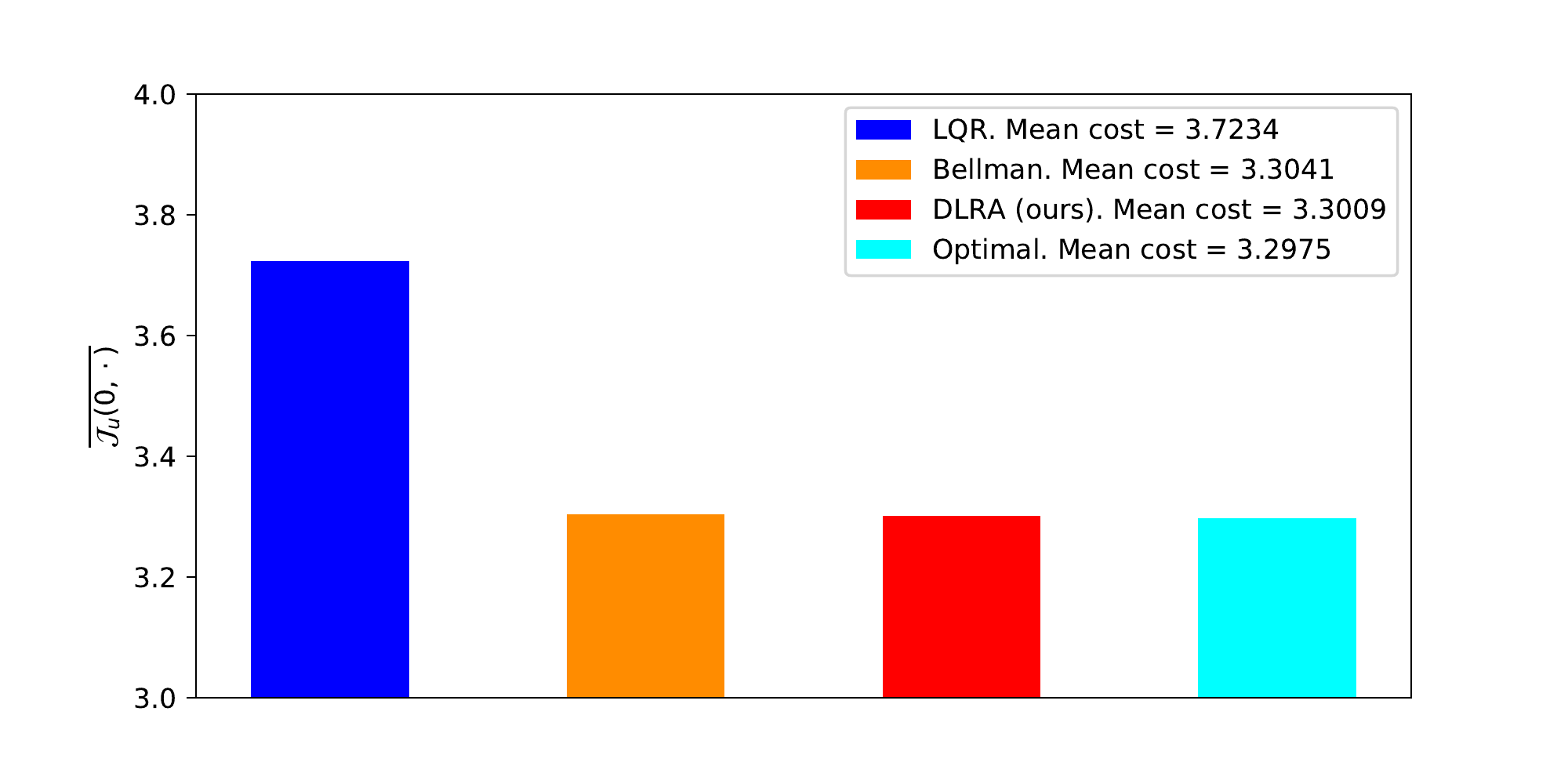}
    \caption{Sample-mean costs $\overline{\mathcal{J}_u(0,\cdot)} \approx \dfrac{1}{N}\sum_{k}^N J_u(0,x^{(k)}_0)$ of uniform initial conditions $x_0^{(k)}$ with the different controllers.}
    \label{fig:costs_const_ic}
\end{figure}

 Figures \ref{fig:control_poly_ic} and \ref{fig:control_const_ic} show the control values $u(t)$ along one trajectory of each type of initial conditions.
 Figures~\ref{fig:costs_poly_ic} and~\ref{fig:costs_const_ic} depict the mean costs over 500 randomly sampled initial conditions in each of the two cases, where we have omitted those initial conditions for which the open-loop solver used to compute the optimal control did not converge.
 Examining the graphs, we note that the Bellman method and the DLRA method achieve similar, almost optimal performance over the chosen test sets.
 Interestingly, the DLRA method actually slightly outperforms the full Bellman method and is often closer to the optimal control trajectories, which for instance can be seen in Figure~\ref{fig:control_const_ic}.
 We attribute this to the generalisation error of the Bellman method: even if the value function approximation should be more accurate -- due to a projection directly onto the manifold -- the associated optimisation is nonlinear and may get stuck in local optima.
 In the DLRA method, we avoid this problem by coping only with linear minimisation problems.
 
\subsection{Computational cost and a hybrid approach}
 
The distinct advantage of the DLRA method is its greatly reduced computational cost.
Table~\ref{tab:comp_costs} contains the computation times for the two methods (Bellman and DLRA), as well as their mean costs on the set of polynomial initial conditions, with the same hyper-parameters and maximal polynomial degrees of $4$, $6$ and $8$, respectively.
We observe that the two methods achieve comparable performance for degrees $6$ and $8$.
However, the DLRA method achieves this performance in roughly one tenth of the time that the Bellman method requires.
We stress again that the version we used is the \textit{fastest} version of the Bellman method available, since we employ the one-step scheme.
As discussed in ~\cite{OS19,OS21}, this method also suffers from error propagation due to a large number of time steps.
Moreover our proposed method projects onto the tangent space, whereas Bellman always tries to project onto the tensor manifold.

\begin{figure}
    \centering
    \includegraphics[scale=0.6]{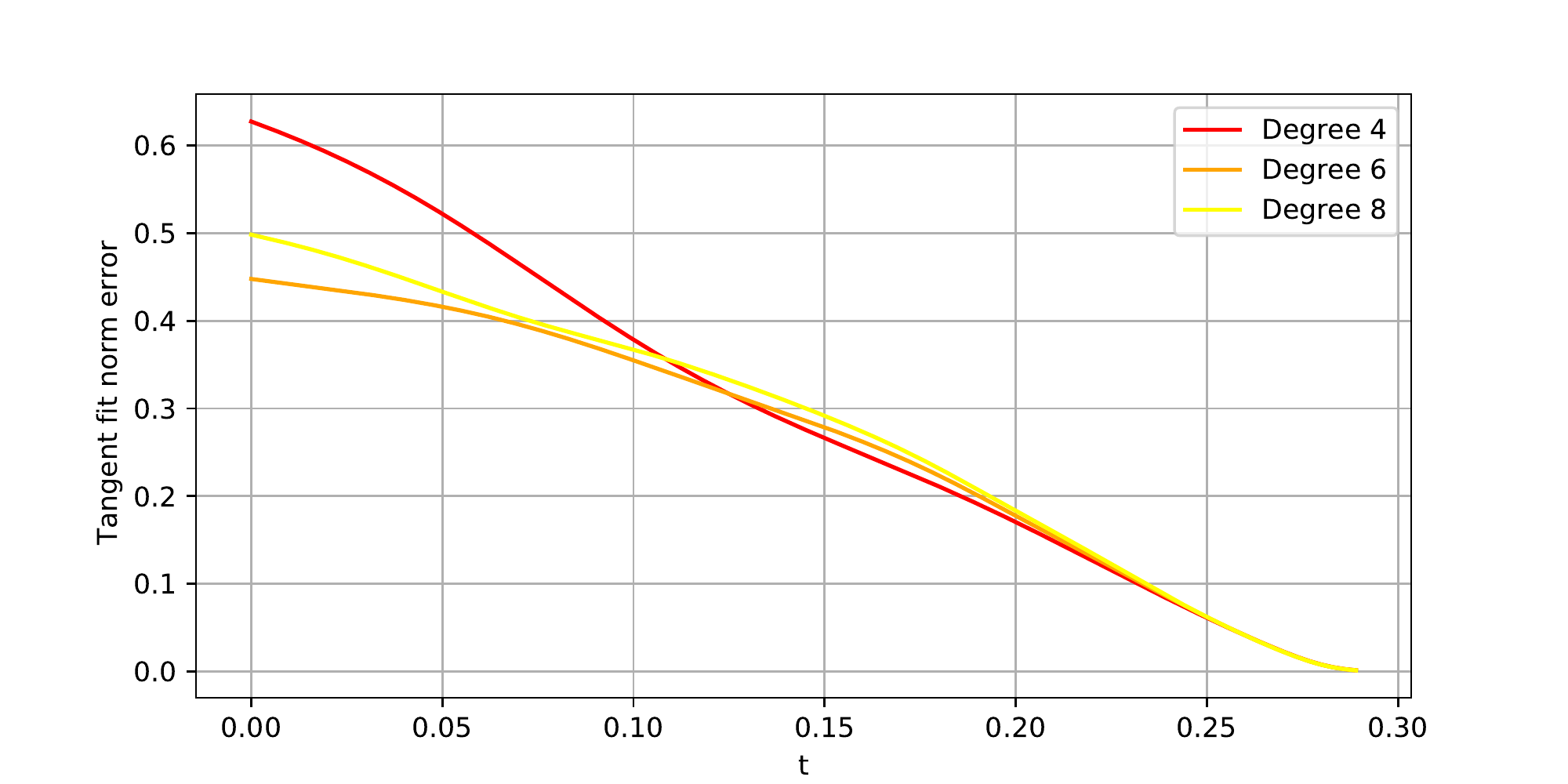}
    \caption{Relative norm error of the last tangent fit (see Appendix \ref{sec:tangent_fit}) in each policy iteration step of the DLRA solver over time and for different polynomial degrees.}
    \label{fig:tangent_res}
\end{figure}

The DLRA method performs significantly worse for a lower polynomial degree of $4$.
We attribute this to an effect that can be seen already for degree $8$ in Figure~\ref{fig:control_const_ic}.
The DLRA controller drifts away from the true optimal control the further it moves away from the terminal time $t=T$.
This error seems to originate from two main factors: for one, the true value function $V^*(t,\cdot)$ successively moves further away from the manifold $\mathcal{M}_{\mathbf{r}}$ even if the terminal condition satisfies $c_T\in\mathcal{M}_{\mathbf{r}}$.
To visualise that the true solution does not stay on the manifold, the relative norm error of the last tangent fit in each policy iteration is plotted over time in Figure~\ref{fig:tangent_res}. Note that these errors should be close to $0$ if the solution to the GHJB equation is an element of the manifold.
Instead, the errors increase monotonically over time.
The second major source of error is the retraction after every Euler step.
Both sources of errors get worse for lower degrees because of the restricted manifold.
Hence, a degree of $4$, which is perfectly feasible for the Bellman method, produces bad results with the DLRA method.
Note that the observed behaviour is expected.
 

\begin{table}[ht]
\centering
\begin{tabular}[t]{lcccccc}
\toprule
 & \multicolumn{2}{c}{Bellman} & \multicolumn{2}{c}{DLRA} & \multicolumn{2}{c}{Hybrid}\\
 & comp. time & mean cost & comp. time & mean cost & comp. time & mean cost \\
\midrule
pol. deg. ~$4$  &3078.44 & 1.8822  & 333.29 & 2.6147 & 909.65
 & 1.8804 \\
pol. deg. ~$6$ & 4270.33 & 1.8801  & 421.52  & 1.8802  & 1851.93
 & 1.8798 \\
pol. deg. ~$8$ & 5967.91 & 1.8800 & 499.96 & 1.8799 & -- & -- \\
\bottomrule
\end{tabular}
\vskip 0.1in
\caption{Computation time of the methods in seconds as well as mean costs of polynomial initial conditions for different maximal polynomial degrees of the basis functions. The mean optimal cost is 1.8793.}
\label{tab:comp_costs}
\end{table}%

This observation leads to a natural formulation of a hybrid method, possibly alleviating the main weaknesses of both methods.
These are the high computational complexity for the Bellman method and error accumulation for the DLRA method.
The hybrid method uses DLRA updates but after each $m$ steps, instead of computing $\hat{V}_i$ with the regular DLRA update, it performs a full Bellman update~\cite{OS21} with an $m$-step scheme
\begin{align}\label{eq:m_step_scheme}
    \hat{V}_i(t_i,x_k) = \int_{t_i}^{t_{i+m}} \ell(t,x_k(t),\alpha(t,x_k(t)) dt  + \hat{V}_{i+m}(t_{{i+m}},x_k(t_{i+m})).
\end{align}

For $m=1$ this method is equivalent to the Bellman method, for $m$ greater than the number of total time steps it is equivalent to the DLRA method.
For any intermediate $m$ it periodically performs one costly but accurate Bellman update in between fast DLRA updates.
Since the maximal number of consecutive DLRA steps is now $m$, the DLRA solver is prevented from drifting too far away from the real solution, before being corrected again by the Bellman update, yielding a new (more accurate) initial condition.
Note in particular that the evaluation of~\eqref{eq:m_step_scheme} does not include any $\hat{V}_j$ computed with the DLRA method.
Hence, after every $m$ steps, the accumulated error of the DLRA steps is reset to $0$.
Globally, only the error of the $m$-step Bellman updates~\eqref{eq:m_step_scheme} accumulates.
 
The results for the hybrid method with $m=10$ are depicted in Table~\ref{tab:comp_costs} for degrees 4 and 6.
We remark that for polynomials of degree  4, the hybrid scheme provides an essential improvement with respect to accuracy when compared to both Bellman and DLRA. 
There is an improvement for degree 6 but compared to DLRA this effect is not pronounced.  Surprisingly, for a sufficiently accurate model, DLRA alone was sufficiently accurate. 
The case of degree 8 is omitted since the DLRA controller is already nearly optimal in that case.
The periodic $10$-step Bellman updates with intermediate DLRA steps are sufficient to outperform the full 1-step Bellman method, but at much lower computational costs.
From the perspective of the DLRA method, the periodic Bellman updates enable the use of more restricted manifolds.

\section{Concluding remarks}
\label{sec:close}

In this paper we present a novel method to approximate optimal feedback laws for optimal control problems.
The proposed method utilizes a tensor train compression to break the curse of dimensionality of a multivariate polynomial ansatz space.
Moreover, it employs an empirical version of the Dirac-Frenkel variational principle to solve the HJB equation.
The method was tested numerically on a canonical benchmark example which is difficult to control with standard methods, and demonstrated to achieve near optimal performance with greatly reduced computation time compared to state-of-the-art methods.

In the experiments it comes as no surprise that the proposed method works quite well for short time intervals.
However, it is striking that we can also observe that with a sufficiently good model -- meaning an adequate polynomial degree in our case -- the method even performs well on a large time horizon. Although we have not considered infinite horizon problems yet, as long as we know stabilizing controls, the present approach probably is applicable as well.
Moreover, for large time horizons we have presented a robust hybrid method.

We would like to point out that the present successful approach strongly exploits the explicit knowledge about the geometry of the considered model class, i.e. (multi-)polynomial tensor trains in our setting.
This advantage is something which cannot be easily transferred to a neural network setting.

We expect the method to also perform favourably with higher dimensional problems, which might be a future research topic.
We predict that this will require some form of rank adaptivity to retain the computational advantage over state-of-the-art methods while achieving similar levels of accuracy.
Rank adaptivity can be incorporated very naturally in the proposed DLRA method: instead of the full retraction onto the manifold $\mathcal{M}_{\mathbf{r}}$ after every step of the solver, one could round the TT based on an adaptive threshold.
Analysing the effect of a changing manifold on the Dirac-Frenkel variational principle might be an interesting topic for future work.

As a second direction, the method could be applied to stochastic optimal control problems.
There, the GHJB equation~\eqref{eq:ghjb} gets an additional Laplacian term $\propto \Delta_x V(t,x)$, turning it into a Kolmogorov-Backward type equation.
Equations of this type for instance govern the time development of observables of It\^o diffusion processes. The application of our method to such problems is currently being investigated.


\section*{Acknowledgements}
Martin Eigel acknowledges the partial support of the DFG SPP 1886 ``Polymorphic Uncertainty Modelling for the Numerical Design of Structures''.
David Sommer acknowledges support by the ProFIT project \textquotedblleft ReLkat – Reinforcement Learning for complex automation engineering\textquotedblright.

\appendix

\section{Details of the empirical risk minimisation}
\label{sec:tangent_fit}

We detail how to reduce the minimisation in~\eqref{eq:emp_risk_min} to a standard system of linear equations.
To achieve this, we use the characterization of the tangent space given by Theorem~\ref{thm:tangent_space} and represent an element $\delta U \in \tau(X)$ of the tangent space as a vector $\mathbf{x}\in\mathbb{R}^{n_X}$.
The first step towards this representation is the parametrisation of the spaces $U_{\mu}^\ell$. 

\subsection{A parametrisation of the tangent space}

By Theorem~\ref{thm:tangent_space}, $U_{\mu}^{\ell}$ is precisely the set of all $r_{\mu-1}n_\mu\times r_\mu$-matrices whose columns are orthogonal to the columns of $L(U_{\mu})$.
Let $QR = L(U_{\mu})$ be the QR decomposition and denote the orthonormal columns of $Q$ by $q_1,\hdots, q_{r_{\mu}}$.
By the Gram-Schmidt procedure, we can expand the columns to an orthonormal basis of $\mathbb{R}^{r_{\mu-1}n_{\mu}}$ and denote the additional vectors by $\hat{q}_1 = q_{r_{\mu}+1},\hdots,\hat{q}_{r_{\mu-1}n_{\mu}-r_{\mu}} = q_{r_{\mu-1}n_{\mu}}$.
Now, let $W_{\mu}\in U_{\mu}^\ell$ and denote its $j$-th column by $w_j$.
Then there are coefficients $c_{j,k}$ such that
\begin{align*}
    w_j = \sum_{k=1}^{r_{\mu-1}n_{\mu}-r_{\mu}}  c_{j,k} \hat{q}_k.
\end{align*}
In total we get $(r_{\mu-1}n_{\mu}-r_{\mu})r_{\mu}$ coefficients $c_{j,k}$, which are stored in a vector
\begin{align*}
    \mathbf{x}_{\mu} &= (c_{1,1},\hdots,c_{1,r_{\mu-1}n_{\mu}-r_{\mu}},c_{2,1},\hdots,c_{2,r_{\mu-1}n_{\mu}-r_{\mu}}, \hdots, \hdots, c_{r_{\mu},1},\hdots, c_{r_{\mu},r_{\mu-1}n_{\mu}-r_{\mu}})^\intercal \\ &\in \mathbb{R}^{(r_{\mu-1}n_{\mu}-r_{\mu})r_\mu} = \mathbb{R}^{r_{\mu-1}n_{\mu} r_\mu-r^2_\mu}.
\end{align*}

From now on we always identify an element of $U_{\mu}^\ell$ with its coefficient vector $\mathbf{x}_{\mu}$.
Elements of $C_d = \mathbb{R}^{r_{d-1}\times n_d\times r_d}$ are represented in the same manner with the only difference that the sum in each column representation goes from $k=1$ to $r_{d-1}n_d$ and the $\hat{q}_i$ can be chosen as the canonical basis in $\mathbb{R}^{r_{d-1}n_d}$.

We eventually can represent an element of the tangent space $\delta U \in \tau(X)$ by the concatenation of its coefficients vectors,
\begin{align*}
    \delta U \cong \mathbf{x} = (\mathbf{x}^{\intercal}_1,\hdots,\mathbf{x}^{\intercal}_d)^\intercal \in \mathbb{R}^{n_X}.
\end{align*}
Since this becomes important when solving the regression problem~\eqref{eq:emp_risk_min} on the tangent space later on, we define a ``lift''
\begin{align*}
    &\mathcal{L}_{\mu}:\mathbb{R}^{r_{\mu-1}n_{\mu} r_{\mu}-r_{\mu}^2}\longrightarrow \mathbb{R}^{r_{\mu-1} n_{\mu} r_{\mu}},
\end{align*}
which maps the coefficient vector of the gauged representation to the vectorised entries of the corresponding tensor in $U_{\mu}^\ell$.
This is achieved by means of a $r_{\mu-1}n_{\mu} r_{\mu}\times(r_{\mu-1}n_{\mu}r_{\mu}-r_{\mu}^2)$-lifting matrix
\begin{align*}
    Z_{\mu} = \begin{pmatrix} 
    Q_{\mu}^\perp & \mathbf{0} & \hdots & \hdots & \mathbf{0} \\ 
    \mathbf{0}  & Q_{\mu}^\perp & \mathbf{0} & & \\ 
    & & \ddots & & \\
    & & & & \mathbf{0}\\
    \mathbf{0} & & & \mathbf{0} & Q_{\mu}^\perp \end{pmatrix}, \qquad Q_{\mu}^\perp = [\hat{q}_1,\hdots, \hat{q}_{r_{\mu-1}n_{\mu} - r_{\mu}}] \in \mathbb{R}^{r_{\mu-1}n_{\mu}\times (r_{\mu-1}n_{\mu}-r_{\mu})}.
\end{align*}
By construction, $Z_{\mu}\mathbf{x}_{\mu}$ is the concatenation of the columns of $L(W_{\mu})$.
Hence, $\mathcal{L}_{\mu}(\mathbf{x}_{\mu})$ can be obtained by $\mathcal{L}_{\mu}(\mathbf{x}_{\mu}) = Z_{\mu} \mathbf{x}_{\mu}$.

\subsection{Solving the system of linear equations}

We examine problem~\eqref{eq:emp_risk_min} in a more general setting.
Let $U\in \mathcal{M}_{\mathbf{r}}$ be $d$-orthogonal and consider the minimisation problem
\begin{align}\label{eq:TT_LGS}
    \min_{T\in\mathcal{T}_U(\mathcal{M}_{\mathbf{r}})} \sum_{k=1}^M \left| T[x^{(k)}] - y^{(k)} \right|^2,
\end{align}
where $(x^{(k)},y^{(k)})_{k=1}^M \subset \mathbb{R}^d\times\mathbb{R}$ is a set of data-target pairs.
This leads to
\begin{align*}
    &~\quad \min_{T\in\mathcal{T}_U(\mathcal{M}_{\mathbf{r}})} \sum_{k=1}^M | T[x^{(k)}] - y^{(k)} |^2 \\
    &= \min_{W\in X} \sum_{k=1}^M | \tau(W)[x^{(k)}] - y^{(k)} |^2 \\
    &= \min_{W\in X} \sum_{k=1}^M \left| \sum_{\mu=1}^d \left[ \sum_{j_1,\hdots,j_d}^{n_1,\hdots,n_d}U_1(j_1)\cdot\hdots\cdot W_{\mu}(j_{\mu}) \cdot\hdots \cdot U_d(j_d) \phi_{j_1}(x_1^{(k)})\hdots \phi_{j_d}(x_d^{(k)})\right] - y^{(k)} \right|^2 \\
    &= \min_{W\in X}  \left\|  \mathfrak{O}(W) - \textbf{y} \right\|_2^2,
\end{align*}
where $\mathbf{y} = (y^{(1)},\hdots,y^{(M)})^\intercal$ and the operator $\mathfrak{O}:X\longrightarrow\mathbb{R}^M$ is defined by
\begin{align}\label{eq:tensor_operator}
    \mathfrak{O}(W) &= \sum_{\mu=1}^d \mathfrak{C}_{\mu}(W_{\mu}), \quad \textnormal{for} \quad W = (W_1,\hdots,W_d),
\end{align}
with
\begin{align*}
    &\mathfrak{C}_{\mu} : \mathbb{R}^{r_{\mu-1}\times n_{\mu} \times r_{\mu}} \longrightarrow \mathbb{R}^M, \\ 
    &(\mathfrak{C}_{\mu}(W_{\mu}))_k = \left[ \sum_{j_1,\hdots,j_d}^{n_1,\hdots,n_d}U_1(j_1)\cdot\hdots\cdot W_{\mu}(j_{\mu}) \cdot\hdots \cdot U_d(j_d) \phi_{j_1}(x_1^{(k)})\hdots \phi_{j_d}(x_d^{(k)})\right].
\end{align*}
Note that $\mathfrak{C}_{\mu}$ is a linear tensor operator in $\mathbb{R}^{M\times r_{\mu-1}\times n_{\mu} \times r_{\mu}}$, which we can transfer into a matrix $C_{\mu} \in \mathbb{R}^{M\times r_{\mu-1}n_{\mu} r_{\mu}}$ by successive unfolding
\begin{align*}
    \mathbb{R}^{M\times r_{\mu-1}\times n_{\mu} \times r_{\mu}} \longrightarrow \mathbb{R}^{M\times r_{\mu-1}n_{\mu} \times r_{\mu}} \longrightarrow \mathbb{R}^{M\times r_{\mu-1}n_{\mu} r_{\mu}}.
\end{align*}
At the first stage, the operator $\mathfrak{C}_{\mu}$ acts on a tensor $W_{\mu}\in U_{\mu}^{\ell}$.
At the second stage, it acts on the left unfolding $L(W_{\mu})$.
And at the third stage, the matrix $C_{\mu}$ acts on the concatenation of the columns of $L(W_{\mu})$.
By the previous section, this concatenation is given by $Z_{\mu} \mathbf{x}_{\mu}$.
We hence have $\mathfrak{C}_{\mu}(W_{\mu}) = C_{\mu} Z_{\mu} \mathbf{x}_{\mu}$, leading to
\begin{align*}
    \mathfrak{O}(W) = \sum_{\mu=1}^d C_{\mu} Z_{\mu} x_{\mu} = A \mathbf{x},
\end{align*}
where $A = \left[ C_1 Z_1,\hdots,C_{d-1}Z_{d-1},C_d \right]$ (note that $Z_d \equiv I_d$).
We have thus transformed~\eqref{eq:TT_LGS} to a standard system of linear equations
\begin{align*}
    \hat{\mathbf{x}} = \argmin_{\mathbf{x}\in\mathbb{R}^{n_X}}\| A\mathbf{x} - \mathbf{y} \|^2 \Longleftrightarrow A^\intercal A \hat{\mathbf{x}} = A^\intercal \mathbf{y}.
\end{align*}
Once a solution is found by standard methods, we recover $W$ from $\hat{\mathbf{x}}$ by reshaping of the component vectors $\mathbf{x}_{\mu}$.

\begin{remark}
We would like to make two remarks about the implementation.
First, note that the matrices $Z_{\mu}$ do not have to be stored in order to compute the product $C_{\mu}Z_{\mu}$ since we can compute
\begin{align*}
    C_{\mu}Z_{\mu} &= [C_{\mu}[0:r_{\mu-1}n_{\mu}]|\hdots | C_{\mu}[(r_{\mu}-1)(r_{\mu-1}n_{\mu}):r_{\mu}r_{\mu-1}n_{\mu}]] \cdot  \textnormal{diag}(Q_{\mu}^\perp,\hdots, Q_{\mu}^\perp)  \\
    &= [C_{\mu}[0:r_{\mu-1}n_{\mu}]Q_{\mu}^\perp|\hdots | C_{\mu}[(r_{\mu}-1)(r_{\mu-1}n_{\mu}):r_{\mu} r_{\mu-1}n_{\mu}]Q_{\mu}^\perp].
\end{align*}

Second, note that $U_{\mu}^{\ell}$ is $0$-dimensional if $r_{\mu-1}n_{\mu}-r_{\mu} = 0$.
In this case, the space consists only of the tensor of constant zeros $\mathbf{0}\in\mathbb{R}^{r_{\mu-1}\times n_{\mu} \times r_{\mu}}$ and hence $W_{\mu} = \mathbf{0}$.
No basis coefficients $\mathbf{x}_{\mu}$ need to be computed.
Consequently, the index $\mu$ can be skipped entirely during optimisation.
By this, $A$ and $\mathbf{x}$ become
\begin{align*}
    A &= [C_1Z_1,\hdots,C_{\mu-1}Z_{\mu-1},C_{\mu+1}Z_{\mu+1},\hdots,C_{d-1}Z_{d-1},C_d], \\
    \mathbf{x} &= (\mathbf{x}^{\intercal}_1,\hdots,\mathbf{x}^{\intercal}_{\mu-1},\mathbf{x}^{\intercal}_{\mu+1},\hdots,\mathbf{x}^{\intercal}_d)^{\intercal}.
\end{align*}
\end{remark}

\printbibliography







\end{document}